\documentclass[10pt]{article} 
\usepackage{amsmath, amssymb, amsthm, mathrsfs}
\usepackage{enumerate}
\usepackage{colordvi}
\usepackage{tikz} 
\usepackage{graphicx}
\usepackage{hyperref}
\newtheorem{claim}{Claim}[section]
\newtheorem{theorem}[claim]{Theorem}
\newtheorem{lemma}[claim]{Lemma}

\newtheorem{corollary}[claim]{Corollary} 
\usepackage{color}

\definecolor{Myred}{cmyk}{0.0,1.0,1.0,0.00}
\definecolor{Mypurple}{rgb}{0.5,0.0,0.5}

\newtheorem{definition}{Definition}

\begin{document}
\title{Spectral convergence of the Laplace operator with Robin boundary conditions on a small hole}
\author
{
Diana Barseghyan$\footnote{corresponding author}$\,\, and Baruch Schneider
}
\date{\small Department of Mathematics, University of Ostrava,  30.dubna 22, Ostrava 70103, Czech Republic\\E-mails:\, diana.barseghyan@osu.cz, baruch.schneider@osu.cz
}
\maketitle

\begin{abstract} In this paper we study a bounded domain with a small hole removed. Our main result concerns the spectrum of the Laplace operator with the Robin conditions imposed at the hole boundary. Moreover we prove that under some suitable assumptions on the parameter in the boundary condition the spectrum of the Laplacian converges in the Hausdorff distance sense to the spectrum of the Laplacian defined on the unperturbed domain.
\end{abstract}
\bigskip

\\


\section{Introduction}
\label{s:intro}

It is a common expectation that small perturbations of the physical situation will lead to only a small change in the spectrum. In the case of domain perturbations, this is largely true for Dirichlet boundary conditions, while the Neumann or Robin case is more delicate. In the recent literature, such questions have already received quite a few answers, starting with the seminal work of Rauch and Taylor
on the spectrum of the Laplace operator of domains with holes \cite{RT75}.  

An excellent shortcut to the recent work on the asymptotic behavior of the eigenvalues of the Laplace operator on the domain with small spherical obstacles imposing the Neumann condition at their boundary and the Dirichlet condition at the rest part of the boundary can be found in  \cite{O83}.

Maz’ya, Nazarov and Plamenewskii, see \cite[Ch.9, vol.I]{MNP00}, have  considered the Laplace operator on the domain with obstacles, imposing the Dirichlet condition on their boundary and have proved the validity of a complete asymptotic expansion for the eigenvalues.

For a survey on more recent research in this subject, we refer the reader to \cite{BN98}, \cite{CP20}, \cite{D03}  where authors have considered the Dirichlet Laplacians on Euclidean domains or manifolds with holes and studied the problems of the resolvent convergence.

The problems with small Neumann obstacles of more general geometry can cause abrupt changes in the spectrum. For example, such an effect is observed when the hole has a "split ring" geometry, see \cite{S15}. The split ring (even if very small) can produce additional eigenvalues that have nothing in common with the eigenvalues of the Neumann Laplacian on the unperturbed domain. The problems with the Neumann obstacles having more general geometry have been studied in \cite{CP20} and later in \cite{BSH22}.

The Robin case for general self-adjoint elliptic operators was considered in \cite{EB20} but with the restriction that the boundary of the unperturbed domain is $C^2$- smooth. In the mentioned work shrinking the hole and scaling properly the parameter in the boundary condition, the authors obtain an operator family that converges, in the norm-resolvent sense, to an operator with a point interaction in the domain without the hole. Results on resolvent convergence for operators with Robin conditions in domains with small holes in higher dimensions were also considered in \cite{BM18}, but in this paper the original domain must again be $C^2$-smooth and the Robin condition was used with the coefficient independent of the hole size.

In this paper we will focus our attention on the spectral properties of a Laplacian defined on a two-dimensional bounded domain with \textit{no additional assumptions} on the smoothness of its boundary with a single hole $K_\varepsilon$ (for a fixed parameter $\varepsilon$) having the Lipschitz boundary.  On the boundary of the original domain we impose the Dirichlet or Neumann condition, and on the boundary of the hole we impose the Robin condition with the coefficient depending on the size of the hole. 

Our main result is the proof of the spectral convergence, in the Hausdorff distance sense, of the spectrum of the Laplacian defined on the perturbed domain to the spectrum of the Laplacian defined on the original domain.

{\it Plan of the paper.} The paper consists of 6 sections, besides this introduction.

 In Section 2 we present the main results, and consider a general theorem, namely Theorem \ref{Second}.  We will use Theorem  \ref{Second} in the proof of Theorem \ref{main} about the spectral convergence for the Laplacians on $\Omega$ and $\Omega\setminus K_\varepsilon$. 

Section 3 contains the main tools of the spectral convergence of operators on varying Hilbert spaces 

In Section 4 and Section 5 we prove our results to which we already alluded.

In Section 6 we give some auxiliary material established in \cite{BSH22}. 

\bigskip

\section{Main results}
\setcounter{equation}{0}

In this section, we present our main results. These results are proven in the following sections.

Let $\Omega\subset \mathbb{R}^2$ be a bounded domain and $K\subset \Omega$ be a compact simply connected set with Lipschitz boundary. We denote 
$\Omega_K:= \Omega\setminus K$. By using Lemma \ref{trace}, the quadratic forms 
\begin{eqnarray}\label{quadratic form}
q_1(u)= \int_{\Omega_K} |\nabla u|^2\,d x\,d y+ \gamma_K\int_{\partial K}|u|_K^2\,d \mu,\quad u\in \mathcal{H}^1(\Omega_K)\,,
\\\label{quadratic form1} 
q_2(u)= \int_{\Omega_K} |\nabla u|^2\,d x\,d y+ \gamma_K\int_{\partial K}|u|_K^2\,d \mu,\quad u\in \mathcal{H}^1(\Omega_K)\cap \mathcal{H}_0^1(\Omega)\,,
\end{eqnarray}
where $\mu$ is the measure on $\partial K$ related to the arc length and $\gamma_K$ is a real number, 
are closed and semi-bounded from below and hence define unique
self-adjoint operators $H_{\Omega_K}(\gamma_K)$ and $\widetilde{H}_{\Omega_K}(\gamma_K)$ which act as the Laplacian on their domains. We will study the question of the convergence of the spectrums of the operators $H_{\Omega_K}(\gamma_K),\widetilde{H}_{\Omega_K}(\gamma_K)$ when the hole $K$ converges to a point. 

We start by a rather important results in the following theorem.

\begin{theorem}\label{main}
Let $\Omega$ be an open bounded domain in $\mathbb{R}^2$. Suppose that ${\mathbb B}_\varepsilon\subset \Omega$ is a ball with center at some point $x_0\in\Omega$ and radius $\varepsilon>0$. Suppose that $K= K_\varepsilon \subset \mathbb B_\varepsilon$ be a bounded simply connected compact set with Lipschitz boundary.
Let $H^N_\Omega$ and $H^D_\Omega$ be the Neumann and Dirichlet Laplacians defined on the unperturbed domain $\Omega$ and $H_{\Omega_{K_\varepsilon}}, \widetilde{H}_{\Omega_{K_\varepsilon}}$ be the operators generated by (\ref{quadratic form}) and (\ref{quadratic form1}) on $\Omega_{K_\varepsilon}$ with the coefficient $\gamma_{K_\varepsilon}>0$ satisfying
\begin{equation}
\label{gamma}
\gamma_{K_\varepsilon}= \varepsilon M_\varepsilon,
\end{equation}
where $M_\varepsilon$ explodes to infinity under the condition that $M_\varepsilon= {o}\left(\frac{1}{\varepsilon^{3/2}}\right)$. Then, for sufficiently small $\varepsilon$, there exists $\eta(\varepsilon)>0$ with $\eta(\varepsilon)\to 0$ as $\varepsilon\to 0$, so that the following spectral convergence occurs 
$$\overline{d} \left(\sigma_{\bullet} \left(H_{\Omega_{K_\varepsilon}}\right), \, \sigma_{\bullet} \left(H^N_\Omega\right)\right) \le \eta(\varepsilon),$$
$$\overline{d} \left(\sigma_{\bullet} \left(\widetilde{H}_{\Omega_{K_\varepsilon}}\right), \, \sigma_{\bullet} \left(H^D_\Omega\right)\right) \le \eta(\varepsilon),$$
where $\overline{d}$ is defined in (\ref{Hausdorff2}) and $\sigma_\bullet (\cdot)$ denotes either the entire spectrum, the essential spectrum, or the discrete spectrum. Furthermore, the multiplicity of the discrete spectrum is preserved.
\end{theorem}

\noindent The previous result motivates the following consequences:


\begin{corollary}
Suppose that $H^N_\Omega$ has purely discrete spectrum denoted by $\lambda^N_k(\Omega)$ (repeated according to multiplicity), and let $\lambda^D_k(\Omega)$ be the discrete spectrum of $H^D_\Omega$. Then the infimum of the essential spectrums of $H_{\Omega_{K_\varepsilon}}, \widetilde{H}_{\Omega_{K_\varepsilon}}$ tend to infinity and there exists $\eta_k(\varepsilon)>0$ with $\eta_k(\varepsilon)\to 0$ as $\varepsilon\to 0$ such that
$$|\lambda^N_k(\Omega)- \lambda_k(\Omega_{K_\varepsilon})|\le \eta_k(\varepsilon), 
$$
$$|\lambda^D_k(\Omega)- \beta_k(\Omega_{K_\varepsilon})|\le \eta_k(\varepsilon)$$
for small enough $\varepsilon$. Here, $\lambda_k(\Omega_{K_\varepsilon})$ and $\beta_k(\Omega_{K_\varepsilon})$ denote the discrete spectrum of $H_{\Omega_{K_\varepsilon}}$ and  $\widetilde{H}_{\Omega_{K_\varepsilon}}$ (below the essential spectrum) repeated according to multiplicity.
\end{corollary}

\begin{corollary} The Hausdorff distance between the spectra of $H_{\Omega_{K_\varepsilon}}$ and $H^N_\Omega$ and $\widetilde{H}_{\Omega_{K_\varepsilon}}$ and $ H^D_\Omega$ converges to zero on any compact interval $[0, \Lambda]$. 
\end{corollary}

The proof of Theorem \ref{main} is based on Theorem \ref{POSTtm} and on the following theorem:

\begin{theorem}\label{Second}
Under the assumptions of Theorem \ref{main}
the operators $H^N_\Omega$ and $H_{\Omega_{K_\varepsilon}}$ are $\delta(\varepsilon)$ close of order $2$ with $\delta(\varepsilon)\to 0$ as $\varepsilon\to 0$. The same is true for the operators $H^D_\Omega$ and $\widetilde{H}_{\Omega_{K_\varepsilon}}$.
\end{theorem}
 
\section{Main tool of the spectral convergence of operators on varying Hilbert spaces}
\setcounter{equation}{0}

For the convenience of the reader, this section begins by reviewing some basic facts that ensure spectral convergence for two operators having different domains. For more information we refer the reader to \cite{P06}. 
To a Hilbert space $H$ with inner product $(\cdot, \cdot)$ and norm $\|\cdot\|$ together with a non-negative, unbounded operator $A$ we associate the scale of Hilbert spaces
\[
H_k:=\mathrm{Dom}((A+I)^{k/2}),\quad \|u\|_k:=\|(A+I)^{k/2}u\|,\,\,k\ge 0,
\]
where $I$ is the identity operator.

We think of $(H', A')$ as some perturbation of $(H, A)$ and want to relax the assumption so that the spectral properties are not the same, but still close.

\begin{definition}(see \cite{P06})\label{Post} 
Suppose we have linear operators
\begin{gather*} J: H \longrightarrow H', \quad\quad J_1: H_1\longrightarrow H_1'\\
J': H'\longrightarrow H,\quad\quad J_1': H_1' \longrightarrow H_1. \end{gather*}
Let $\delta> 0$ and $k\ge 1$.
We say that $(H, A)$ and $(H', A')$ are $\delta$-close of order $k$ iff the following conditions are fulfilled:

\begin{gather}\label{1} \|J f- J_1 f\|_0\le \delta \|f\|_1, \\\label{2}  |(J f, u)-(f, J' u)|\le \delta \|f\|_0 \|u\|_0,\\\label{3}  \|u-J J' u\|_0\le \delta \|u\|_1,\\\label{4} \quad \|Jf\|_0\leq 2\|f\|_0,\,\,\, \|J'u\|_0\le 2\|u\|_0,
\\\label{5'} \|(f- J'J f)\|_0\le \delta \|f\|_1,\\
\label{6} \| J' u- J_1'u\|_0\le \delta \|u\|_1\\\label{7}  |a(f, J_1' u)- a'(J_1 f, u)|\le \delta \|f\|_k \|u\|_1,\end{gather}
for all $f, u$ in the appropriate spaces. Here, $a$ and $a'$ denote the sesquilinear forms associated to $A$ and $A'$.
\end{definition}

We denote by $d_{\mathrm{Haussdorff}}(A, B)$ the Hausdorff distance for subsets $A, B\subset\mathbb{R}$:
\begin{equation}\label{Haus.}
d_{\mathrm{Haussdorff}}(A, B):=\mathrm{max} \left\{\underset{a\in A} {\mathrm{sup}}\,d(a, B),\,\underset{b \in B} {\mathrm{sup}}\, d(b, A)\right\},
\end{equation}
where $d(a, B):= \mathrm{inf}_{b\in B} |a - b|$. We set
\begin{equation}\label{Hausdorff2}\overline{d} (A, B):= d_{\mathrm{Hausdorff}}\left((A+I)^{-1}, (B+I)^{-1}\right)\end{equation}
for closed subsets of $[0, \infty)$. 

For the next result, which comes from the work of O. Post \cite{P06}, we have the following spectral convergence theorem in terms of the distance $\overline{d}$.
\begin{theorem}\cite{P06}\label{POSTtm}
There exists $\eta(\delta)>0$ with $\eta(\delta)\to 0$ as $\delta\to 0$ such that
\begin{equation}\label{Post}\overline{d}(\sigma_\bullet(A), \, \sigma_\bullet (A'))\le \eta(\delta)\end{equation}
for all pairs of non-negative operators and Hilbert spaces $(H, A)$ and $(H', A')$ that are $\delta$-close. Here, $\sigma_\bullet (A)$ denotes either the entire spectrum, the essential or the discrete spectrum of $A$. Also, the multiplicity of the discrete spectrum, 
$\sigma_{\mathrm{disc}}$, is preserved, i.e. if $\lambda\in \sigma_{\mathrm{disc}}(A)$ has multiplicity $\mu>0$, then there exist $\mu$ eigenvalues (not necessarily all distinct) of the operator $A'$ belonging to the interval $(\lambda-\eta(\delta), \lambda+\eta(\delta))$.
\end{theorem}

We now turn to the proof of Theorem \ref{Second}. Since the proof is almost the same for both the Dirichlet and Neumann cases, we will restrict ourselves to the Neumann case.
The only difference is Lemma \ref{Delta'}, but its validity for the Dirichlet case can be easily checked from \cite{BSH22}.

\section{Proof of Theorem\,\ref{Second}}

At this stage, we divide the proof into two steps. 
\bigskip
\\
{\bf Step \,1.}\,{\it Construction of the mappings $J, J', J_1, J_1'$.}
\bigskip
\\
We can apply the technique of \cite{BSH22}. It is easy to see that
$H= L^2(\Omega),\,H'= L^2(\Omega_{K_\varepsilon}),\,A= A'=-\Delta$,
$H_1,\,H_1'$ correspond to Sobolev spaces $\mathcal{H}^1(\Omega)$ and $\mathcal{H}^1(\Omega_{K_\varepsilon})$ and
$H_2= \mathrm{Dom}(H^N_\Omega)$. The norm $\|\cdot\|_0$ corresponds with the
$L^2$ norm,
\begin{eqnarray*}&&\|u\|_1=\left(\|u\|^2_0+\|\nabla u\|_0^2+\gamma_{K_\varepsilon}\int_{\partial K_\varepsilon}|u|^2\,d\mu\right)^{1/2},\quad\text{for}\quad u\in\mathcal{H}^1(\Omega_{K_\varepsilon}),\\&&\|u\|_1=\left(\|u\|^2_0+\|\nabla u\|_0^2\right)^{1/2},\quad\text{for}\quad u\in\mathcal{H}^1(\Omega)\\&&\text{and}\quad
\|f\|_2=\|-\Delta f+f\|_0.\end{eqnarray*} 

We define $J u=J_1 u=u|_{\chi_{\Omega_{K_\varepsilon}}}$ for all $u\in H$ and  
{$J' u =u \chi_{\Omega_{K_\varepsilon}}$ for all $u\in H'$.}

Now let us construct the mapping $J_1': H_1' \to H_1$.
Without loss of generality, assume that the ball $\mathbb B_\varepsilon$ mentioned in Theorem \ref{main} and Theorem \ref{Second} is centered at the origin. Let $\epsilon\in (\varepsilon, 2\varepsilon)$ be a number to be chosen later and let $\mathbb B_\epsilon \supset \mathbb B_\varepsilon$ be the ball again centered on the origin and radius 
$\epsilon$, $\Omega_\epsilon:=\Omega\backslash \mathbb B_\epsilon$.

We will first construct the mapping $J_1'$ first for smooth functions. For each $v\in  C^\infty(\Omega_{K_\varepsilon})$  we define
\begin{equation*}J_1' v:=\begin{cases} v,\quad \text{on} \quad \Omega_\epsilon, \\ \frac{r}{\epsilon} \tilde{v}(\epsilon, \varphi), \quad \text{on} \quad \mathbb B_\epsilon,\end{cases}\end{equation*}
where $\tilde{v}(r, \varphi)=v(r \cos \varphi, r \sin \varphi)$.

Now let us construct the mapping $J_1' u$ for any $u\in H_1'$.  Using the approximation method described in 
\cite[Thm.2, 5.3.2]{E10}, for the fixed sequence $\{\eta_k\}_{k=1}^\infty$ converging to zero we construct the sequence 
$v_{\eta_k}\in C^\infty(\Omega_{K_\varepsilon})$ which satisfies
\begin{equation}\label{converging1}
\int_{\Omega_{K_\varepsilon}}|\nabla (u-v_{\eta_k})|^2\,d xd y+\int_{\Omega_{K_\varepsilon}}|u-v_{\eta_k}|^2\,d xd y< \eta_k\,\|u\|_1^2.\end{equation}

To deal with $\int_{\partial K_\varepsilon}|u|^2\,d\mu$ we will use the trace inequality \cite{E10}. We present it immediately after an auxiliary result on Lipschitz bounds \cite{E10}.

\begin{lemma}\label{mu}
Let $\omega$ be a bounded two-dimensional open set with a Lipschitz boundary. Then there exist $\delta>0$ and $\mu\in C^\infty(\overline{\omega})$ such that
$$\mu\cdot \nu\ge \delta, \quad\text{a.e. on}\quad \partial \omega,$$
where $\nu$ is the normal to $\partial \omega$.
\end{lemma}

\begin{lemma}\label{trace}
Let $v\in \mathcal{H}^1(\omega)$, where $\omega$ is a two-dimensional domain with Lipschitz boundary.
Then there exists a constant $K>0$ depending on $\omega$ such that
$$\int_{\partial \omega} 
|v|^2d\mu\le K \int_\omega\left(\delta|\nabla v|^2+\frac{1}{\delta}|v|^2\right)\,d x\,d y,$$
where $K$ depends only on the norm of $\mu$ in $C^1(\overline{\Omega})$ and $\delta\in (0, 1)$.
\end{lemma}

Combining the above lemma with $\delta=\frac{1}{2}$ and the inequality (\ref{converging1}) we get
\begin{eqnarray*}&&\int_{\partial K_\varepsilon}|u-v_{\eta_k}|^2\,d \mu\\&&\le 2K\,\left(\int_{\Omega_{K_\varepsilon}}|\nabla (u-v_{\eta_k})|^2\,d x\,d y+\int_{\Omega_{K_\varepsilon}}|u-v_{\eta_k}|^2\,d x\,d y\right)\\&&<2K  \eta_k \|u\|_1^2.\end{eqnarray*}
Therefore
\begin{equation}
\label{converging}
\|u-v_{\eta_k}\|_1<\sqrt{(1+2K \gamma_{K_\varepsilon})\,\eta_k}\,\|u\|_1.
\end{equation} 

Let us mention that in view of the inequalities (\ref{J1'est.}) and (\ref{J 1''}), which will be proved later, and the construction of the function $J_1'$, it follows that for any smooth function $v$ the integrals
$\int_\Omega |\nabla J_1'v|^2\,d x\,d y$ and $\int_\Omega |J_1'v|^2d x\,d y$ can be estimated from above 
by $\|v\|_1^2$ multiplied by some constant. Combining this with the Lemma \ref{trace} we get
\begin{equation}\label{inequality} \|J_1'v\|_1^2\le \overline{C}(\varepsilon) \|v\|_1^2,\end{equation}
where $\overline{C}(\varepsilon)$ is some constant.

Due to the positivity of the coefficient $\gamma_{K_\varepsilon}$ and Lemma \ref{trace} the completeness of the space $H_1$ is equivalent to the completeness of the Sobolev space $\mathcal{H}^1$.
Thus, using (\ref{converging}) and (\ref{inequality}) we can define
\begin{equation}\label{approximation}
J_1'u= \underset{k\to \infty}{\mathrm{lim}} \,J_1'v_{\eta_k}.
\end{equation}
\\
{\bf Step\, 2.}\, \it {The conditions (\ref{1})-(\ref{7}) hold for the mappings $J, J', J_1, J_1'$}.

\bigskip
\rm  Indeed, we have that the estimates (\ref{1})-(\ref{4}) are satisfied with $\delta=0$. 
\\

We now prove (\ref{5'}), i.e., {\it under the assumptions stated in Theorem \ref{Second}, inequality (\ref{5'}) is satisfied with $\delta=\mathcal{O}\left(\sqrt{\frac{\varepsilon}{\gamma_{K_\varepsilon}}+\varepsilon}\right)$ for small enough $\varepsilon$. Thus, given (\ref{gamma}), $\delta$ converges to zero as $\varepsilon\to 0$.}

Given our construction, we have 
\begin{eqnarray}\nonumber
\|f- J' J f\|_0^2=\int_\Omega |f- J' f|_{\Omega_{K_\varepsilon}}|^2\,d x\,d y=\int_\Omega |f- f\chi_{\Omega_{K_\varepsilon}}|^2\,d x\,d y\\\label{start 6}=\int_{K_\varepsilon}|f|^2\,d xd y\le \int_{\mathbb{B}_\varepsilon}|f|^2\,d x\,d y.
\end{eqnarray}

To complete the proof of (\ref{5'}), we use the following lemma, applied with $\eta=\varepsilon$ and
$\Gamma=\emptyset$:
\begin{lemma}\label{auxiliary}
Let $\Omega$ be an open bounded domain in $\mathbb{R}^2$. Suppose that $\mathbb{B}_\epsilon\subset \Omega$ is a ball with center at some point $x_0\in\Omega$ and radius $\epsilon>0$. Suppose that $\Gamma\subset \mathbb{B}_\epsilon$ be a bounded simply connected compact set with Lipschitz boundary.
Then for any function $u\in\mathcal{H}^1(\mathbb{B}_\epsilon\setminus \Gamma)$ the following inequality holds
$$\int_{\mathbb{B}_\epsilon\setminus \Gamma}|u(x, y)|^2\,d xd y\le C_1\left(\frac{\epsilon}{\gamma_\Gamma}+\epsilon\right)\|u\|_1^2,$$
where $C_1>0$ is a constant that depends on the distance between the boundary of 
$\mathbb{B}_\epsilon$ and the boundary of $\Omega$. 
\end{lemma}

Let us pass to (\ref{6}), i.e. {\it under the assumptions given in Theorem\,\ref{Second} inequality (\ref{6}) is satisfied with $\delta=\mathcal{O}\left(\sqrt{\frac{\varepsilon}{\gamma_{K_\varepsilon}}+\varepsilon}\right)$ for small enough $\varepsilon$. Thus, in view of (\ref{gamma}), $\delta$ converges to zero as $\varepsilon\to 0$.}

Using the construction of $J'$ and $J_1'$ and the completeness of $C^\infty(\Omega_{K_\varepsilon})$ in space $H_1'$ it is sufficient to prove (\ref{6}) for $u\in C^\infty(\Omega_{K_\varepsilon})$.

Considering that $J' u=0$ on $K_\varepsilon$, $J'u= u$ on $\Omega_{K_\varepsilon}$ and $J_1' u= u$ on $\Omega\setminus \mathbb{B}_\epsilon$ one has
\begin{gather}
\nonumber
\|J' u- J_1' u\|_0^2= \int_{\Omega_{K_\varepsilon}}|u-J_1'u|^2\,d xd y+ \int_{K_\varepsilon}|J_1' u|^2\,d xd y\\\nonumber =\int_{\mathbb{B}_\epsilon\setminus K_\varepsilon}|u-J_1'u|^2\,d xd y+ \int_{K_\varepsilon}|J_1' u|^2\,d xd y
\\\label{ineq}\le 2\int_{\mathbb{B}_\epsilon\setminus K_\varepsilon}|u|^2\,d xd y+2\int_{\mathbb{B}_\epsilon}|J_1' u|^2\,d xd y+ \int_{K_\epsilon}|J_1' u|^2\,d xd y.
\end{gather}

To estimate the first integral on the right-hand side of (\ref{ineq}), we use the Lemma \ref{auxiliary} with 
$\eta=\epsilon$ and $\Gamma=K_\varepsilon$.

Now let us examine the second term. Passing to polar coordinates, we get
\begin{eqnarray*}
\nonumber &&\int_{\mathbb B_\epsilon}|J_1'u|^2\,d x\,d y = \int_0^\epsilon \int_0^{2\pi} r \left|\frac{r}{\epsilon}
\tilde{u}(\epsilon, \varphi)\right|^2\,d r\,d \varphi\\&&\le\int_0^\epsilon \int_0^{2\pi} r |\tilde{u}(\epsilon, \varphi)|^2\,d r\, d\varphi\le \epsilon \int_0^\epsilon \int_0^{2\pi} 
|\tilde{u}(\epsilon, \varphi)|^2\,d r\,d \varphi\\&&= \epsilon^2 \int_0^{2\pi} 
|\tilde{u}(\epsilon, \varphi)|^2\,d \varphi.
\end{eqnarray*}
Taking into account that $\epsilon \int_0^{2\pi} 
|\tilde{u}(\epsilon, \varphi)|^2\, d \varphi$ coincides with the curvilinear integral 
$\int_{\partial \mathbb{B}_\epsilon} |u|^2\, d \mu$,  the above bound leads to
\begin{equation}\label{ugamma}
\int_{\mathbb B_\epsilon}|J_1'u|^2\,d xd y\le\epsilon
\int_{\partial \mathbb{B}_\epsilon}|u|^2\,d \mu.
\end{equation}

In view of Lemma \ref{trace} applied with $\delta=\frac{1}{2}$, one estimates from above the right-hand side of (\ref{ugamma}) by $2K \epsilon\int_{\Omega\setminus \mathbb{B}_\epsilon}(|\nabla u|^2+|u|^2)\,d x\,d y $. Then use the following obvious bound which holds due to the positivity of $\gamma_{K_\varepsilon}$:
\begin{equation}\label{H1}\int_{\Omega\setminus \mathbb{B}_\epsilon}(|\nabla g|^2+|g|^2)\,d x\,d y\le \|g\|_1^2,\quad\text{for all}\quad g\in \mathcal{H}^1(\Omega_{K_\epsilon}),\end{equation}
and the fact that $\epsilon\le 2\varepsilon$ we have
\begin{equation}\label{J1'est.}
\int_{\mathbb B_\epsilon}|J_1'u|^2\,d xd y\le 4K\varepsilon \|u\|_1^2.
\end{equation}
Indeed, we have
$$\int_{K_\varepsilon}|J_1'u|^2\,d x\,d y\le 4K\varepsilon \|u\|_1^2,$$
since
$$\int_{K_\varepsilon}|J_1'u|^2\,d x\,d y\le\int_{B_\epsilon}|J_1'u|^2\,d x\,d y.$$

Combining the above inequality together with (\ref{J1'est.}), the right-hand side of (\ref{ineq}) can be estimated as follows
$$
\|J' u- J_1'u\|_0^2\le\left(2C_1\left(\frac{\varepsilon}{\gamma_{K_\varepsilon}}+\varepsilon\right)+12K \varepsilon\right) \|u\|_1^2 
$$
which proves (\ref{6}) with $\delta=\mathcal{O}\left(\sqrt{\frac{\varepsilon}{\gamma_{K_\varepsilon}}+\varepsilon}\right)$ for small enough $\varepsilon$.\\

We now give the proof of the estimate (\ref{7}), i.e. {\it under the assumptions given in Theorem \ref{Second}, the inequality (\ref{7}) holds with $k=2$ and $\delta=\mathcal{O}\left(\varepsilon^{1/6}+ \gamma_{K_\varepsilon}^{1/2} \varepsilon^{1/4}\right)$ for sufficiently small 
$\varepsilon$. Thus, in view of (\ref{gamma}), $\delta$ converges to zero as $\varepsilon\to 0$.}\\

As before without loss of generality suppose that $u\in C^\infty(\Omega_{K_\varepsilon})$.
Since $J_1'u=u$ on $\Omega_{\mathbb{B}_\epsilon}$ and $J_1u=u$ on $\Omega_{K_\varepsilon}$ we have
\begin{eqnarray}\nonumber&&
|a(f, J_1'u)- a'(J_1f, u)|\\\nonumber&&=\left|\int_\Omega\nabla f\, \overline{\nabla (J_1'u)}\,d x\,d y-\int_{\Omega_{K_\varepsilon}}\nabla (J_1f)\,\overline{\nabla u}\,d x\,d y
- \gamma_{K_\varepsilon} \int_{\partial K_\varepsilon}(J_1f)\, \overline{u}\,d \mu\right|\\\nonumber&&
=\biggl|\int_{\Omega_{\mathbb B_\epsilon}}\nabla f\,\overline{\nabla u}\,d x\,d y+\int_{\mathbb B_\epsilon}\nabla f\,\overline{\nabla (J_1'u)}\,d x\,d y-\int_{\Omega_{K_\varepsilon}}\nabla (J_1f)\,\overline{\nabla u}\,d x\,d y- \gamma_{K_\varepsilon} \int_{\partial K_\varepsilon}(J_1f)\, \overline{u}\,d \mu\biggr|\\\nonumber&&=\biggl|\int_{\mathbb B_\epsilon}\nabla f\,\overline{\nabla (J_1'u)}\,d x\,d y-\int_{\mathbb B_\epsilon\backslash K_\varepsilon}
\nabla f\, \overline{\nabla u}\,d x\,d y- \gamma_{K_\varepsilon} \int_{\partial K_\varepsilon}(J_1f)\, \overline{u}\,d \mu\biggr|\\\label{7''}&&\le\left|\int_{\mathbb B_\epsilon}\nabla f\,\overline{\nabla (J_1'u)}\,d x\,d y\right|+\left|\int_{\mathbb B_\epsilon\backslash K_\varepsilon}
\nabla f\,\overline{\nabla u}\,d x\,d y\right|+ \gamma_{K_\varepsilon} \left|\int_{\partial K_\varepsilon}(J_1f)\, \overline{u}\,d \mu\right|.\end{eqnarray}

Let us estimate each term of (\ref{7''}). Starting with the first term we get
\begin{equation}\label{8}
\left|\int_{\mathbb B_\epsilon}\nabla f\,\overline{\nabla (J_1'u)}\,d xd y\right|\le\left(\int_{\mathbb B_\epsilon}|\nabla f|^2\,d xd y\right)^{1/2}\left(\int_{\mathbb B_\epsilon}|\nabla (J_1'u)|^2\,d xd y\right)^{1/2}.
\end{equation}
 
Since $f\in \mathcal{H}^2_{\mathrm{loc}}(\Omega)$, then using Lemma \ref{magnetic} (see Appendix) applied with domain $\Omega'$ such that $\overline{\Omega'}\subset \Omega$ and $\mathbb{B}_\epsilon\subset \Omega'$, and the fact that $\epsilon\le 2\varepsilon$,
the first term on the right hand side of (\ref{8}) can be estimated as follows\begin{equation}\label{rhs9}\int_{\mathbb B_\epsilon}|\nabla f|^2\,dx\,dy\le 2^{4/3}C_2\varepsilon^{4/3}\int_{\Omega'}|-\Delta f+f|^2\,d x\,d y\le2^{4/3}C_2\varepsilon^{4/3}\|f\|_2^2.
\end{equation}

To proceed further with the proof of an upper bound of (\ref{8}), we need to estimate the integral $\int_{\mathbb B_\epsilon}|\nabla (J_1'u)|^2\,d x\,d y$. Passing to polar coordinates we get
\begin{eqnarray}\nonumber&&\int_{\mathbb B_\epsilon} |\nabla J_1' u|^2\,d x\,d y=
\int_{\mathbb B_\epsilon} \left(\left|\frac{\partial (J_1' u)}{\partial x}\right|^2+ \left|\frac{\partial (J_1' u)}{\partial y}\right|^2 \right)\,d x\,d y\\\nonumber&&
= \int_0^\epsilon \int_0^{2\pi}
r \left|\frac{1}{\epsilon} \tilde{u}(\epsilon, \varphi)\,\cos \varphi -\frac{1}{\epsilon} \frac{\partial \tilde{u}}{\partial \varphi}(\epsilon, \varphi)\, \sin \varphi\right|^2\,d r\,d \varphi\\\nonumber&&+ \int_0^\varepsilon \int_0^{2\pi}
r \left|\frac{1}{\epsilon} \tilde{u}(\epsilon, \varphi)\,\sin \varphi +\frac{1}{\epsilon} \frac{\partial \tilde{u}}{\partial \varphi}(\epsilon, \varphi)\, \cos \varphi\right|^2\,d r\,d \varphi\\\nonumber&& \le \frac{4}{\epsilon} \int_0^\epsilon \int_0^{2\pi} \left(|\tilde{u}(\epsilon, \varphi)|^2+\left|\frac{\partial \tilde{u}}{\partial \varphi}(\epsilon, \varphi)\right|^2\right)\,d r\,d \varphi\\\label{J 1'}&& =4 \int_0^{2\pi}|\tilde{u}(\epsilon, \varphi)|^2d\varphi+ 4 \int_0^{2\pi}\left|\frac{\partial \tilde{u}}{\partial \varphi}(\epsilon, \varphi)\right|^2\,d \varphi.
\end{eqnarray}

As in the proof of (\ref{6}), we find that
$$\int_0^{2\pi}|\tilde{u}(\epsilon, \varphi)|^2\,d \varphi=\frac{1}{\epsilon}\,\int_{\partial \mathbb B_\epsilon} 
|u|^2\,d \mu.$$

Thus, using Lemma \ref{trace} applied with $\delta=\frac{1}{2}$ and inequality (\ref{H1}), we have
\begin{equation}\label{first term}
\int_0^{2\pi}|\tilde{u}(\epsilon, \varphi)|^2\,d \varphi\le \frac{2K}{\epsilon} \|u\|_1^2.
\end{equation}

Now we come to the second term of (\ref{J 1'}). Given the Lemma \ref{auxiliary1} (see Appendix) used for $g=u$, there exists a number $\tau\in (\varepsilon, 2\varepsilon)$ such that 
$$
\int_0^{2\pi} \left|\frac{\partial \tilde{u}}{\partial \varphi}(\tau, \varphi)\right|^2\,d \varphi \le 4 \int_{\mathbb B_{2\varepsilon}\setminus \mathbb B_\varepsilon}(|\nabla u|^2+|u|^2)\,d xd y,
$$
where $\tilde{u}(r, \varphi):= u(r\cos\varphi, r\sin\varphi)$ and $\tau\in (\varepsilon, 2\varepsilon)$ is some number.

If $\tau$ belongs to the interval 
$(\varepsilon, 3\varepsilon/2]$, then we take $\epsilon$ as the supremum of all such numbers in
$(\varepsilon, 3\varepsilon/2]$. In the opposite case if $\tau\in (3\varepsilon/2, 2\varepsilon)$, then let $\epsilon$ be the infimum of such numbers. Since $u$ is a smooth function, the above inequality is satisfied with $\tau=\epsilon$. 

Combining this together with (\ref{H1}) we get
$$\int_0^{2\pi}\left|\frac{\partial \tilde{u}}{\partial \varphi}(\epsilon, \varphi)\right|^2\,d \varphi\le4 \|u\|_1^2.$$
 
Thus, by virtue of (\ref{J 1'}), (\ref{first term}), the above estimate and using the fact that $\epsilon\ge \varepsilon$, we have
\begin{equation}\label{J 1''}\int_{\mathbb B_\epsilon} |\nabla J_1' u|^2\,d x\,d y\le\left(\frac{8K}{\varepsilon}+ 16\right)\|u\|_1^2.\end{equation}

Finally, using the above bound and inequality (\ref{rhs9}), we estimate the right-hand side of (\ref{8}) as follows
\begin{equation}\label{f.term}
\left|\int_{\mathbb B_\epsilon}\nabla f\,\overline{\nabla (J_1'u)}\,d x\,d y\right|\le 2^{13/6}(C_2(K+2\varepsilon))^{1/2}\varepsilon^{1/6} \|f\|_2\|u\|_1. 
\end{equation}

Let us now consider the second term in (\ref{7''}). By virtue of  (\ref{H1}) and (\ref{rhs9}) we get
\begin{eqnarray}\nonumber\left|\int_{\mathbb B_\epsilon\backslash K_\varepsilon}
\nabla f\,\overline{\nabla u}\,dx\,dy\right|\le\left(\int_{\mathbb B_\epsilon}
|\nabla f|^2\,dx\,dy\right)^{1/2}\left(\int_{\mathbb B_\epsilon\backslash K_\varepsilon}|\nabla u|^2\,d x\,d y\right)^{1/2}\\\label{second term}\le 2^{2/3}(C_2)^{1/2}\varepsilon^{2/3}\|f\|_2\|u\|_1.
\end{eqnarray} 

Finally we move on to the third term in (\ref{7''}). We have
\begin{eqnarray}\nonumber
\left|\int_{\partial K_\varepsilon}(J_1f)\, \overline{u}\,d \mu\right|\le \left(\int_{\partial K_\varepsilon}|J_1f|^2\,d\mu\right)^{1/2} \left(\int_{\partial K_\varepsilon} |u|^2\,d \mu\right)^{1/2}\\\label{partial.}= \left(\int_{\partial K_\varepsilon}|f|^2\,d \mu\right)^{1/2} \left(\int_{\partial K_\varepsilon} |u|^2\,d \mu\right)^{1/2}.
\end{eqnarray}

Let us first find the appropriate estimate for the first integral on in the right hand side of (\ref{partial.}).

Let $\Pi(d)=(-d, d)^2,\,d>0$, be the maximum square belonging to $\Omega$ and containing $K_\varepsilon$.
For almost all $x_0$ belonging to the projection of $K_\varepsilon$ on the axis $X$, let  $y(x_0)\in (-d, d)$ be the point such that 
$$
|f(x_0, y(x_0))|^2\le \frac{1}{d} \int_{-d}^d|f(x_0, t)|^2\,d t\,.
$$

Let us fix any $y\in (-d, d)$. Without loss of generality suppose that $y> y(x_0)$. Then  
\begin{gather}\nonumber 
|f(x_0, y)|^2=\left|f(x_0, y(x_0))+ \int_{y(x_0)}^y\frac{\partial f}{\partial t}(x_0, t)\,d t\right|^2\\\nonumber\le 2\left|f(x_0, y(x_0))\right|^2+ 2\left|\int_{y(x_0)}^y\frac{\partial f}{\partial t}(x_0, t)\,d t\right|^2
\\\nonumber\le\frac{2}{d} \int_{-d}^d |f(x_0, t)|^2\,d t+ 2(y-y(x_0))\int_{y_0(x_0)}^y\left|\frac{\partial f}{\partial t}(x_0, t)\right|^2\,d t\\\label{K1}\le\frac{2}{d} \int_{-d}^d|f(x_0, t)|^2\,d t+ 4d\int_{-d}^d\left|\frac{\partial f}{\partial t}(x_0, t)\right|^2\,d t\,.
\end{gather}

Without loss of generality assume that the boundary of the unperturbed set which is $\frac{1}{\varepsilon} K_\varepsilon$ is parametrized as $(x, y_1(x)),\,x\in(-1, 1)$, where $y_1$ is some $C^1$-smooth function. Then the parameterization  of the boundary of $K_\varepsilon$ coincides with $(x, \varepsilon y_1(x/\varepsilon)),\,x\in(-\varepsilon, \varepsilon)$.  

Integrating $|f(x, y)|^2$ over $\partial K_\varepsilon$ and using the inequality (\ref{K1}) we get
\begin{gather}\nonumber
\int_{\partial K_\varepsilon}|f(x, y)|^2\,d \mu=\int_{-\varepsilon}^\varepsilon |f(x, \varepsilon y_1(x/\varepsilon))|^2\,\left(1+y_1'^2(x/\varepsilon)\right)^{1/2}\,d x
\\\nonumber\le\left(1+ \|y_1'\|_{L^\infty(-1, 1)}^2\right)^{1/2} \int_{-\varepsilon}^\varepsilon\biggl(\frac{2}{d} \int_{-d}^d|f(x, t)|^2\,d t+ 4d\int_{-d}^d\left|\frac{\partial f}{\partial t}(x, t)\right|^2\,dt\biggr)\,d x\\\label{C}\le 2\mathrm{max}\left\{1/d, 2d\right\} \left(1+ \|y_1'\|_{L^\infty(-1, 1)}^2\right)^{1/2}\int_{\Omega_\varepsilon} (|f|^2+ |\nabla f|^2)\,d xd y\,,
\end{gather}
where $\Omega_\varepsilon:=(-\varepsilon, \varepsilon)\times (-d, d)$. 

To proceed with a proof we need the following auxiliary result \cite{MK06}:

\begin{lemma}\label{marchenko}
Let $\Pi'\subset \mathbb{R}^n$ be a convex set and let $G$ and $Q$ be arbitrary measurable sets in $\Pi'$ with $\mu\, (G)\ne 0$.
Then, for all $v\in \mathcal{H}^1(\Pi')$, the following inequality holds:
\begin{eqnarray}\label{Marchenko'}
&&\int_Q |v|^2\,d x\,d y\\ \nonumber
&&\le \frac{2\mu\, (Q)}{\mu \,(G)} \int_G |v|^2\,d x\,d y
+\frac{C(n) (d(\Pi'))^{n+1} (\mu\, (Q))^{1/n}}{\mu \,(G)} \int_{\Pi'} |\nabla v|^2\, d x\,d y,
\end{eqnarray}  
where $d(\Pi')$ is the parameter of $\Pi'$, $\mu$ is the Lebesque measure on $\mathbb{R}^n$, and the constant $C(n)$ depends only on the dimension of $\mathbb{R}^n$.
\end{lemma}

Let $G=\Pi'=\Pi(d)$ and $Q=\Omega_\varepsilon$. Then using the Lemma \ref{marchenko} for the function $f$ we get
\begin{eqnarray*}&&\int_{\Omega_\varepsilon}|f|^2\,d x\,d y\\&&\le \frac{8 \varepsilon d}{\mu(\Pi(d))}\int_{\Pi(d)}|f|^2\,d x\,d y+\frac{2 C(2) (d(\Omega))^3 \sqrt{\varepsilon d}}{\mu(\Pi(d))}\int_{\Pi(d)}|\nabla f|^2\,d x\,d y,
\\
&&\int_{\Omega_\varepsilon}|\nabla f|^2\,d x\,d y\\&&\le \frac{8 \varepsilon d}{\mu(\Pi(d))}\int_{\Pi(d)}|\nabla f|^2\,d x\,d y+
\frac{2 C(2) (d(\Omega))^3 \sqrt{\varepsilon d}}{\mu(\Pi(d))}\int_{\Pi(d)}\left(\left|\nabla\left(\frac{\partial f}{\partial x}\right)\right|^2+ \left|\nabla\left(\frac{\partial f}{\partial y}\right)\right|^2\right)\,d x\,d y\\\nonumber
&&\le \frac{8 \varepsilon d}{\mu(\Pi(d))}\int_{\Pi(d)}
|\nabla f|^2\,d x\,d y+ \frac{2 C(2) (d(\Omega))^3 \sqrt{\varepsilon d}}{\mu(\Pi(d))}\int_{\Pi(d)}\left(\left|\frac{\partial^2 f}{\partial x^2}\right|^2+ \left|\frac{\partial^2 f}{\partial y^2}\right|^2+ 2\left|\frac{\partial^2 f}{\partial x\partial y}\right|^2\right)\,d x\,d y\,.\end{eqnarray*}

With the above bounds, one can show that sufficiently small values of $\varepsilon$, the following is true 
\begin{equation}\label{Pi}
\int_{\Omega_\varepsilon}(|\nabla f|^2+ |f|^2)\,d xd y\le \tilde{c}\,\sqrt{\varepsilon}\,\|f\|_{\mathcal{H}^2(\Pi(d))},
\end{equation}
where 
$$
\tilde{c}= \frac{2 C(2) (d(\Omega))^3 \sqrt{d}}{\mu(\Pi(d))}=\frac{C(2)(d(\Omega))^3}{2d \sqrt{d}}
$$
and $\|f\|_{\mathcal{H}^2(\Pi(d))}$ means the Sobolev $\mathcal{H}^2(\Pi(d))$ norm of $f$.

Next we need the following interior regularity theorem \cite{B16}:
\begin{theorem}(Interior Regularity Theorem.)
Suppose that $h\in \mathcal{H}^1(\Omega)$ is a weak solution of 
$-\Delta h=w$. Then $h\in \mathcal{H}^2_{\mathrm{loc}}(\Omega)$ and for each $\Omega_0\subset  \Omega$ there exists a constant $c= c(\Omega_0)$ independent of $h$ and $w$ such that: 
\begin{equation}\label{Elliptic}\|h\|_{\mathcal{H}^2(\Omega_0)}\le c \left(\|h\|_{L^2(\Omega)}+\|w\|_{L^2(\Omega)}\right).
\end{equation}
\end{theorem}

In view of the above theorem the right-hand side of (\ref{Pi}) can be estimated as follows
\begin{equation}\label{Pi'}\int_{\Omega_\varepsilon}(|\nabla f|^2+ |f|^2)\,d x\,d y\le c\,\tilde{c}\,\sqrt{\varepsilon}\,\int_\Omega (|f|^2+|\Delta f|^2)\,d x\,d y,\end{equation}
with some constant $c=c(d)$ does not depend on $\varepsilon$.

By virtue of inequalities (\ref{C}), (\ref{Pi'}) and Lemma \ref{Delta'} (see Appendix) we have
\begin{eqnarray*}\int_{\partial K_\varepsilon}|f(x, y)|^2\,d \mu
\le 2c\,\tilde{c}\, \mathrm{max}\left\{1/d, 2d\right\} \left(1+ \|y_1'\|_{L^\infty(-1, 1)}^2\right)^{1/2}\, \sqrt{\varepsilon}\,\,\int_\Omega (|f|^2+|\Delta f|^2)\,d xd y\\\le
2c\,\tilde{c}\, \mathrm{max}\left\{1/d, 2d\right\} \left(1+ \|y_1'\|_{L^\infty(-1, 1)}^2\right)^{1/2}\, \sqrt{\varepsilon}\,\|f\|_2^2\,.
\end{eqnarray*}

The above combined the fact that 
$$\gamma_{K_\varepsilon} \int_{\partial K_\varepsilon} |u|^2\,d\mu\le \|u\|_1^2$$ estimates the right-hand side of the inequality (\ref{partial.}) as follows 
\begin{eqnarray*}&&\left|\int_{\partial K_\varepsilon}(J_1f)\, \overline{u}\,d \mu\right|\\&&\le \left(\frac{2c\,\tilde{c}\, \mathrm{max}\left\{1/d, 2d\right\}}{\gamma_{K_\varepsilon}}\right)^{1/2} \left(1+ \|y_1'\|_{L^\infty(-1, 1)}^2\right)^{1/4}\,\varepsilon^{1/4}\,\|f\|_2 \|u\|_1\,.
\end{eqnarray*}

By virtue of (\ref{f.term}), (\ref{second term}) and the above inequality the right-hand side of inequality (\ref{7''}) satisfies
$$
\mathrm{r.h.s.}(\ref{7''})=\mathcal{O}\left(\varepsilon^{1/6}+ \gamma_{K_\varepsilon}^{1/2} \varepsilon^{1/4}\right)\,\|f\|_2\, \|u\|_1,\quad \varepsilon\to 0,
$$
which ends the proof.\qed

\section{Proof of Lemma \ref{auxiliary}}
\setcounter{equation}{0}

Let $M_\epsilon$ be the subset of $\mathbb{B}_\epsilon\setminus\Gamma$ such that for every point $(x_0, y_0)\in M_\epsilon$ the line $l_{x_0}:=\{x=x_0\}$ intersects the boundary of $\Gamma$ at least once. Let $(x_0, y_1(x_0))$ be the point of intersection of $l_{x_0}$ with $\partial \Gamma$ and let $(x_0, y_2(x_0))\in\Omega$ be the point such that the open interval connecting $(x_0, y_1(x_0, y_0))$ and $(x_0, y_2(x_0, y_0))$ belongs to $\Omega$ and has an empty intersection with $\Gamma$. 
Without loss of generality suppose that $y_1(x_0)< y_2(x_0)$. Then for any $u\in\mathcal{H}^1(\mathbb{B}_\varepsilon\setminus \Gamma)$ and almost all $(x_0, y_0)\in M_\epsilon$ we have
$$u(x_0, y_0)= u(x_0, y_1(x_0, y_0))+ \int_{y_1(x_0, y_0)}^y \frac{\partial u}{\partial t}(x_0, t)\,d t.$$

Let $P_{\epsilon}$ denote the projection of $M_\epsilon$ onto the axis $\mathbb{X}$.  We get
\begin{eqnarray}\nonumber
&&\int_{M_\epsilon}|u(x, y)|^2\,d x\,d y\\\nonumber&&\le 2\int_{M_\epsilon}|u(x, y_1(x, y))|^2\,d x\,d y+ 2\int_{M_\epsilon}\left|\int_{y_1(x, y)}^y \frac{\partial u}{\partial t}(x, t)\,d t\right|^2\,d x\,d y\\&&\nonumber= 2\int_{P_{\epsilon}} \int_{\{y:\, (x, y)\in  M_\epsilon\}}|u(x, y_1(x))|^2\,d xd y\\\nonumber&&+ 2\int_{P_{\epsilon}} \int_{\{y:\,(x, y)\in M_\epsilon\}}\left|\int_{y_1(x, y)}^y \frac{\partial u}{\partial t}(x, t)\,d t\right|^2\,d x\,d y\\
&&\nonumber\le 2\int_{P_{\epsilon}} \int_{\{y:\, (x, y)\in  M_\epsilon\}}|u(x, y_1(x))|^2\,d x\,d y\\\nonumber&&+ 2\int_{P_{\epsilon}} \int_{\{y:\,(x, y)\in M_\epsilon\}}(y-y_1(x))\int_{y_1(x)}^y \left|\frac{\partial u}{\partial t}(x, t)\right|^2\,d t\,d x\,d y\\\nonumber&&\le 4\epsilon
\int_{P_{\epsilon}}|u(x, y_1(x))|^2\,d x+ 8\epsilon^2\int_{P_{\epsilon}} \int_{y_1(x)}^{y_2(x)}\left|\frac{\partial u}{\partial t}(x, t)\right|^2\,d x\,d t\\\nonumber&&\le 4\epsilon
\int_{P_{\epsilon}}|u(x, y_1(x))|^2\,d x+ 8\epsilon^2\int_{\mathbb{B}_\epsilon\setminus\Gamma}\left|\frac{\partial u}{\partial t}(x, t)\right|^2\,d x\,d t\\\label{M}&&\le4\epsilon
\int_{P_{\epsilon}}|u(x, y_1(x))|^2\,d x+ 8\epsilon^2 \|u\|_1^2.\end{eqnarray} 

Since the set $\{(x, y_1(x))\},\,x\in P_{\epsilon}$, is the part of $\Gamma$ then
\begin{eqnarray*}\int_{P_\epsilon}|u(x, y_1(x))|^2\,d x\le \int_{P_{\epsilon}} \left(1+y_1'(x)^2\right)^{1/2}|u(x, y_1(x))|^2\,d x\\= \int_{\Gamma\cap\{(x, y): x\in P_{\epsilon}\}}|u|^2\,d \mu\le \int_\Gamma |u|^2\,d \mu\,.
\end{eqnarray*}

Returning to the inequality (\ref{M}) and combining the above bound together with the fact that
$$\int_\Gamma |u|^2\,d \mu\le \frac{1}{\gamma_\Gamma}\|u\|_1^2$$
we get
\begin{equation}\label{first result}
\int_{M_\epsilon}|u(x, y)|^2\,d xd y\le 4\left(\frac{\epsilon}{\gamma_\Gamma}+ 2\epsilon^2\right)\|u\|_1^2\,.
\end{equation}
Now let us go to the subset $\left(\mathbb{B}_\epsilon\setminus\Gamma\right)\setminus M_\epsilon$.
For any $(x_0, y_0)\in\left(\mathbb{B}_\epsilon\setminus\Gamma\right)\setminus M_\epsilon$ let $(x_0, y_3(x_0))$ with $y_3(x_0)<y_0$, be a point of intersection of line $l_{x_0}$ with the boundary of $\Omega$. One can easily check that there is
$y_4(x_0)\in (y_3(x_0), y_0)$ such that
$$|u(x_0, y_4)|\le\frac{1}{\sqrt{y_0-y_3(x_0)}}\sqrt{\int_{y_3(x_0)}^{y_0}|u(x_0, z)|^2\,d z}.$$

Therefore 
\begin{eqnarray}\nonumber
|u(x_0, y_0)|^2= \left|u(x_0, y_4(x_0))+\int_{y_4(x_0)}^{y_0}\frac{\partial u}{\partial z}(x_0, z)\,d z\right|^2\\\nonumber\le 2|u(x_0, y_4(x_0))|^2+ 2(y_0-y_4(x_0))\int_{y_4(x)}^y\left|\frac{\partial u}{\partial z}(x, z)\right|^2\,d y\\\nonumber\le \frac{2}{y_0-y_3(x_0)}\int_{y_3(x_0)}^{y_0}|u(x_0, z)|^2\,d z
+ 2(y_0-y_3(x_0))\int_{y_3(x_0)}^{y_0}\left|\frac{\partial u}{\partial z}(x_0, z)\right|^2\,d z\\\nonumber\le 
\frac{2}{\mathrm{dist}((x_0, y_0), \partial \Omega)}\int_{y_3(x_0)}^{y_0}|u(x_0, z)|^2\,d z
+ 2 \mathrm{diam}(\Omega)\int_{y_3(x_0)}^{y_0}\left|\frac{\partial u}{\partial z}(x_0, z)\right|^2\,d z\\\nonumber\le \widetilde{C} \left(\int_{y_3(x_0)}^{y_0}|u(x_0, z)|^2\,d z
+ \int_{y_3(x_0)}^{y_0}\left|\frac{\partial u}{\partial z}(x_0, z)\right|^2\,d z\right)\,,
\end{eqnarray} 
where  $\mathrm{diam}(\Omega)$ is the diameter of $\Omega$ and $\mathrm{dist}((x_0, y_0), \partial\Omega)$ is the distance between $(x_0, y_0)$ and the boundary of $\Omega$ and $$\widetilde{C}=2 \mathrm{max}\left\{\frac{1}{\mathrm{dist}((x_0, y_0), \partial \Omega)},\,  \mathrm{diam}(\Omega)\right\}.$$

Let $P_\epsilon'$ be the projection of $\left(\mathbb{B}_\epsilon\backslash \Gamma\right)\setminus M_\epsilon$ onto axis $\mathbb{X}$. So we get
\begin{eqnarray*}
\int_{\left(\mathbb{B}_\epsilon\backslash \Gamma\right)\setminus M_\epsilon}|u(x, y)|^2\,d xd y\le \widetilde{C} \int_{\left(\mathbb{B}_\epsilon\backslash \Gamma\right)\backslash M_\epsilon} \int_{y_3(x)}^y\left(\left|\frac{\partial u}{\partial z}(x, z)\right|^2+|u(x, z)|^2\right)\,d zd xd y\\
= \widetilde{C}\int_{P'_\epsilon}\int_{\{y: (x, y)\in \left(\mathbb{B}_\epsilon\backslash \Gamma\right)\backslash M_\epsilon\}} \int_{y_3(x)}^y\left(\left|\frac{\partial u}{\partial z}(x, z)\right|^2+|u(x, z)|^2\right)\,d zd xd y\,.\end{eqnarray*} 

Then finally
$$
\int_{\left(\mathbb{B}_\epsilon\backslash \Gamma\right)\setminus M_\epsilon}|u(x, y)|^2\,d xd y\le 2\widetilde{C} \epsilon  \int_{\Omega\backslash \Gamma}(|\nabla u|^2+ |u|^2)\,d xd y\le 2\widetilde{C} \epsilon \|u\|_1^2\,.
$$

This together (\ref{first result}) proves the lemma with $$C_1= 
\mathrm{max}\{2\widetilde{C}+ 8\varepsilon, 4\}< \mathrm{max}\{2\widetilde{C}+ 8, 4\}.$$
\qed

\section{Appendix}
\setcounter{equation}{0}

In this section we mention several useful lemmas proved in \cite{BSH22}.

\begin{lemma}\label{magnetic}
Let $\Omega'$ be an open bounded domain in $\mathbb{R}^2$ and let $p\in \Omega'$ be some fixed point. Suppose that $\mathbb{B}_\epsilon\subset \Omega'$ is a ball with center at $p$ and radius $\epsilon> 0$.
For any function $g\in \mathcal{H}^2(\Omega')$ the following estimate takes place
\[
\int_{\mathbb B_\epsilon}|\nabla g|^2\,d xd y\le C_2\epsilon^{4/3}\int_{\Omega'}|-\Delta g+g|^2\,d xd y,
\]
with the constant $C_2$ depending on $\Omega$.
 \end{lemma}

\bigskip

\begin{lemma}\label{auxiliary1}
Let $\mathbb B_{2\varepsilon}$ and $\mathbb B_\varepsilon$ be the balls centered on the origin and the radii $\varepsilon$ and $2\varepsilon$. Let $g\in \mathcal{H}^1(\mathbb B_{2\varepsilon}\setminus \mathbb B_\varepsilon)$. Then there exists $\tau\in (\varepsilon, 2\varepsilon)$ such that
$$\int_0^{2\pi} \left|\frac{\partial \tilde{g}}{\partial \varphi}(\tau, \varphi)\right|^2\,d \varphi \le 4 \int_{\mathbb B_{2\varepsilon}\setminus \mathbb B_\varepsilon}(|\nabla g|^2+|g|^2)\,d xd y,$$
where $\tilde{g}(r, \varphi):= g(r\cos\varphi, r\sin\varphi)$.
\end{lemma}

\bigskip

\begin{lemma}\label{Delta'}
For any $z\in \mathrm{Dom}(H_\Omega^N)$ the following estimate is valid
$$\int_\Omega |-\Delta z+ z|^2\,d xd y \ge \int_\Omega (|\Delta z|^2+ |z|^2)\,d xd y.$$
\end{lemma}

\subsection*{Acknowledgements}
D.B. acknowledges support from the Czech Science Foundation (GACR), the project 21-07129S.


\bigskip


\begin{thebibliography}{10}

\bibitem{B16} L.~ Beck, Elliptic Regularity Theory, Springer International Publishing, Switzerland, 2016

\bibitem{BSH22} D.~Barseghyan, B.~Schneider and L.H.~Hai, Neumann Laplacian in a perturbed domain, The Mediterranean Journal of Mathematics (2022), 19(126), 1--17. 

\bibitem{BM18} D.I.~ Borisov,  A.I.~ Mukhametrakhimova, On norm resolvent convergence for elliptic operators in multi-dimensional domains with small holes, Journal of Mathematical Sciences 232 (2018), 283--298.

\bibitem{BN98} M.~Balzano and L.~Notarantonio, On the asymptotic behavior of Dirichlet problems in a Riemannian manifold less small random holes, Rendiconti del Seminario Matematico della Universita di Padova 100 (1998), 249--282.

\bibitem{CP20} A.~Colette, O.~Post, Wildly perturbed manifolds: norm resolvent and spectral convergence, Journal of Spectral Theory, European Mathematical Society, In press.

\bibitem{D03} D.~Daners, Dirichlet problems on varying domains, Journal of Differential Equations 188(2003), 591--624.

\bibitem{E10} L. C.~Evans, Partial differential equations, American Mathematical Society, 2010.

\bibitem{EB20} P.~Exner, D.~Borisov, Approximation of point interactions by geometric
perturbations in two-dimensional domains, Bulletin of Mathematical Sciences(2022) 2250003.

\bibitem{MK06} V. A.~Marchenko, E. Ya.~Khruslov, Homogenization of partial differential equations, Progress in Mathematical Physics, 46, Birkh\"auser, Boston, 2006.

\bibitem{MNP00}  V. G.~Maz’ya, S. A.~Nazarov, B. A.~ Plamenewskii, Asymptotic Theory of Elliptic Boundary Value Problems in Singularly Perturbed Domains, I, II. Operator Theory: Advances and Applications, vols. 111, 112. Birkh\"{a}user,Basel (2000) (translation of the original in German published by Akademie Verlag 1991)

\bibitem{O83}  S.~Ozawa,  Spectra of domains with small spherical Neumann boundary, Proceedings of the Japan Academy, Series A, Mathematical Sciences 58(5)  (1982), 190--192.
 
\bibitem{P06} O.~Post, Spectral convergence of quasi-one-dimensional spaces, Annales Henri Poincare 7 (2006), 933-- 973.

\bibitem{RT75} J.~Rauch,  M.~Taylor, Potential and scattering theory on wildly perturbed domains, Journal of Functional Analysis 18(1975), 27--59.

\bibitem{S15} B.~Schweizer, The low-frequency spectrum of small Helmholtz resonators,
Proceedings of the Royal Society A: Mathematical, Physical and Engineering Sciences: 20140339, 2015.


\end{thebibliography}
\end{document}